\documentclass[10pt]{amsart}

\usepackage{color,graphics,graphicx,shortvrb}
\usepackage{epsfig}
\usepackage{newlfont}
\usepackage{caption}

\usepackage{amssymb,amsmath,latexsym}
\usepackage{mathptmx}
\usepackage{subfig}
\usepackage{bibentry}
\usepackage{dsfont}

\newtheorem{thm}{Theorem}[section]

\newtheorem{rem}{Remark}[section]

\newcommand{\R}{\mathbb{R}}
\newcommand{\N}{\mathbb{N}}

\newcommand{\la}{\lambda}
\newcommand{\vf}{\varphi}

\newcommand{\p}{\partial}
\begin{document}

\title[ Pointwise stabilization of Cascade ODE-PDE]{Exponential Stabilization of Cascade ODE-Reaction-Diffusion PDE by Pointwise Actuation}

\author{Habib Ayadi}
\address{Universit\'e de Kairouan. Institut Sup\'erieur des Math\'ematiques Appliqu\'ees et de l'Informatique, Avenue Assad Iben Fourat - 3100 Kairouan, Tunisia}
\email{hayadi26@laposte.net}
\begin{abstract}
In this paper, we are concerned with the state feedback stabilization  of ODE-PDE cascade systems governed by a linear ordinary differential
equation  and the $1-d$ reaction-diffusion equation  posed on a bounded interval. In contrast to the previous works in the literature where the control acts
at the boundary, the control for the entire system acts at an inside point of the PDE domain whereas the PDE acts in the linear ODE by a
Neumann connection. We use the infinite dimensional backstepping design to convert system under consideration to an exponentially target system. By invertibility
of the design and Lyapunov analysis, we prove the well posedness and exponential stability of such system.
\end{abstract}
\maketitle

\tableofcontents

\thispagestyle{empty}

\section{Introduction}
Originally developed for finite dimensional control systems governed by ODE \cite{Kris95},  the first extension of the backstepping method appeared
in \cite{Coron98} and \cite{Kris00} for parabolic PDE.  Later, in \cite{Liu03} and \cite{Kris04}, the authors have introduced an invertible integral
transformation that transforms the original parabolic PDE into an asymptotically stable one. Recently, the
backstepping method is used to design a feedback control law for coupled PDE-ODE (see the textbook \cite{Krisbook} and references therein).
Many problems of state and output feedback stabilization for  coupled  ODE-Heat has been solved \cite{Tang11}, \cite{Tang11SCL},
\cite{Kris09heat} and ODE- Wave \cite{Kris10JFI}, \cite{Zhoua12},  to cite few. In all those works, the actuator acts at the left or the right
boundary of the PDE domain by Dirichlet or Neumann actuation. In contrast to the previous works in the literature, in this paper, we propose a model where
the controller acts at a point inside the domain  of the PDE subsystem. More precisely, we consider an ODE-PDE cascade system  governed by
a linear ordinary differential equation  and  the $1-d$ reaction-diffusion equation  posed on a bounded interval $(0,l)$ where the control of the entire system
is located at a point $\xi\in (0,l)$  under a transmission conditions and the PDE acts  in the  ODE by a Neumann connection. We prove that for all $\xi\in (0,l)$,
system is well posed and exponentially stable in the sense of the $H^1$-norm.
In recent years, by using semi semigroups theory, a lot of researches have been devoted to the study of distributed plants with pointwise actuator.
Surprisingly, the control properties of those systems are very different depending of the location of the actuator and the type of boundary conditions.
(see \cite{Am01hetu} for the wave equation and \cite{Am00tu} for the beam equation, to cite few.)
In the other hand, by using the backstepping method,  the exponential stability for reaction-diffusion equation \cite{Woi13} is proved
if $\frac{l}{\xi}=\frac{p}{q}$ co-prime where $p$ odd . However, in \cite{Fa17}, a feedback law  has been proposed to achieve exponential stability
for wave equation for all $\xi\in(0,l)$.
The paper is organized as follows. In section 2, the problem is stated and the main result of this paper
is summarized in Theorem \ref{thm1} . In section 3,  the backstepping method is used to derive the  state-feedback control law.
Section 4 is devoted to the proof of the main Theorem \ref{thm1}. In section 5, we present the conclusion and the future work.

\section{Problem Formulation and Main Result}
Let $l>0$, $\xi\in (0,l)$ and $\lambda>0$, we consider the following cascade ODE-PDE system
\begin{equation}\label{sys1}
\left\lbrace\begin{array}{lll}
\dot{X}(t)=AX(t)+B\frac{\p u}{\p x}(0,t),\quad t>0,\\
\frac{\p u}{\p t}(x,t) = \frac{\partial^2 u}{\partial x^2}(x,t)+\lambda u(x,t) ,\quad t>0, \quad x\in(0,\xi)\cup (\xi,l),\\
u(\xi^-,t) = u(\xi^+,t) ,\quad t>0,\\
\frac{\partial u}{\partial x}(\xi^-,t)-\frac{\partial u}{\partial x}(\xi^+,t)=U(t) ,\quad t>0,\\
X(0)= X^0,\\
u(0,t)  = u(l,t)= 0,\quad t>0,\\
u(x,0)=u^0(x) ,\quad x\in (0,l).
\end{array}\right.
\end{equation}
where $X(t)\in \R^n$ is the state of the ODE subsystem, $u(x, t)\in \mathbb R$ is the state of the reaction-diffusion subsystem,
$U(t)\in \R$ is the control input to the entire system acting in the interior point $x=\xi$ of the PDE domain  $(0, l)$,
and  $A\in \mathbb{R}^{n\times n}$ , $B \in \mathbb{R}^{n\times 1}$  such that the pair $(A, B)$ is  controllable.
When $U(t)=0$, the PDE subsystem is equivalent to the following system
\begin{equation}\label{open}
\left\lbrace\begin{array}{lll}
\frac{\p u}{\p t}(x,t) = \frac{\partial^2 u}{\partial x^2}(x,t)+\lambda u(x,t) ,\quad t>0, \quad x\in(0,l)\\
u(0,t)  = u(l,t)= 0,\quad t>0,\\
u(x,0)=u^0(x) ,\quad x\in (0,l).
\end{array}\right.
\end{equation}
System (\ref{open}) is unstable with arbitrarily many unstable eigenvalues for large $\lambda$. Thus the open-loop system (\ref{sys1}) is unstable
for large $\lambda$. The control objective is to exponentially stabilize  system (\ref{sys1}) around its zero equilibrium.
Dividing  the domain $[0,l]$ into  the subdomains $[0, \xi]$ and $[\xi, L]$, we can reformulated system (\ref{sys1}) as  the following cascade
ODE-transmission system
\begin{equation}\label{sys2}
\left\lbrace\begin{array}{lll}
\dot{X}(t)=AX(t)+B\frac{\p u_1}{\p x}(0,t),\quad t>0,\\
\frac{\p u_1}{\partial t}(x,t) = \frac{\p^2 u_1}{\partial x^2}(x,t)+\lambda u_1(x,t) ,\quad t>0, \quad x\in(0,\xi),\\
\frac{\p u_2}{\p t}(x,t) = \frac{\p^2 u_2}{\partial x^2}(x,t)+\lambda u_2(x,t) ,\quad t>0, \quad x\in(\xi,l),\\
 u_1(\xi,t) = u_2(\xi,t),\\
 \frac{\p u_1}{\p x}(\xi,t)-\frac{\p u_2}{\p x }(\xi,t)=U(t) ,\quad t>0,\\
 u_1(0,t)=u_2(l,t)=0, \quad t>0,\\
 X(0)= X^0,\\
 u_1(x,0)=u_1^0(x) ,\quad x\in (0,\xi),\\
 u_2(x,0)=u_2^0(x) ,\quad x\in (\xi,l).
 \end{array}\right.
\end{equation}
where $$u(x,t)=
\left\lbrace\begin{array}{lll}
u_1(x,t), \quad x\in [0,\xi], \quad t\geq 0,\\
u_2(x,t), \quad x\in [\xi,l], \quad t\geq 0.
\end{array}\right.$$
The backstepping method is to use the transformations
\begin{equation}\label{w1}
w_1(x,t)= u_1(x,t)-\int_0^x k_1(x,y)u_1(y,t)dy+\vf(x)X(t), \quad x\in [0,\xi], \quad t\geq 0,
\end{equation}
and
\begin{equation}\label{w2}
w_2(x,t)= u_2(x,t)+\int_x^l k_2(x,y)u_2(y,t)dy, \quad x\in [\xi,l], \quad t\geq 0,
\end{equation}
such that with the feedback law
\begin{eqnarray}\label{fedb}
U(t)&=&\Big(k_1(\xi,\xi)-k_2(\xi,\xi)\Big)u_1(\xi,t)+\int_0^{\xi}\frac{\p k_1}{\partial x}(\xi,t)u_1(y,t)dy\nonumber\\
&& +\int_{\xi}^1 \frac{\p k_2}{\p x}(\xi,t) u_2(y,t)dy-\vf^{'}(\xi)X(t),
\end{eqnarray}
where the gain kernels\, $k_1(x,y)\in\R$ and \,$k_2(x,y)\in\R$ and the gain function $\vf(x)^T\in\R^n$ are appropriately
chosen to transform system (\ref{sys2}) into the following exponentially stable  target system

\begin{equation}\label{target}
\left\lbrace\begin{array}{lll}
\dot{X}(t)=(A+BK)X(t)+B\frac{\p w_1}{\p x}(0,t),\quad t>0,\\
\frac{\partial w_1}{\partial t}(x,t) = \frac{\partial^2 w_1}{\partial x^2}(x,t),\quad t>0, \quad x\in(0,\xi),\\
\frac{\partial w_2}{\partial t}(x,t) = \frac{\partial^2 w_2}{\partial x^2}(x,t),\quad t>0, \quad x\in(\xi,l),\\
 w_1(\xi,t) = w_2(\xi,t),\\
 \frac{\partial w_1}{\partial x}(\xi,t)=\frac{\partial w_2}{\partial x}(\xi,t) ,\quad t>0,\\
 w_1(0,t)=w_2(l,t)=0, \quad t>0,\\
 X(0)=X^0,\\
 w_1(x,0)=w_1^0(x) ,\quad x\in (0,\xi),\\
 w_2(x,0)=w_2^0(x) ,\quad x\in (\xi,l),
 \end{array}\right.
\end{equation}
where $K\in\R^{1\times n}$ is such that the matrix $A+BK$ is Hurwitz.
Once the transformations (\ref{w1}) and (\ref{w2}) (namely $k_1(z,y)$ , $k_2(z,y)$ \, and $\vf(x)^T$) are found, we use their invertibility
and exponential stability of (\ref{target}) to get that of the original plant (\ref{sys1}) with the feedback law (\ref{fedb}) acting at $\xi$.
Along this paper, the Euclidean norm of a vector $X$ in $\R^n$ and the $L^2$-norm of
a function $u$ in $L^{2}(a,b)$ are denoted by $|X|$ and
$$\|u\|=\Big(\int_a^b u^2(x)dx\Big)^{\frac{1}{2}},$$ respectively.  Let $ H=\mathbb{R}^n\times L^2(0,\xi)\times L^2(\xi,l) $ be the state space of
the  system (\ref{sys2}). It is obvious that the vector space $H$ equipped  with  its norm
\begin{equation}\label{normH}
\|(X,u_1,u_2)\|_H=\Big(|X|^2+\|u_1\|^2+\|u_2\|^2\Big)^{\frac{1}{2}},
\end{equation}
is a Hilbert space. As far as that goes, we denote by  $Y= \R^n\times H_L^1(0,\xi)\times H_R^1(\xi,l)$ the  dense subspace of $H$
endowed with the norm
$$\|(X,u_1,u_2)\|_Y=\Big(|X|^2+\|u_1\|_{H^1(0,\xi)}^2+\|u_2\|_{H^1(\xi,l)}^2\Big)^{\frac{1}{2}},$$
where $H_L^1(0,\xi):=\{v\in H^1(0,\xi), v(0)=0\}$ \,\,and\,\, $H_R^1(\xi,l):=\{v\in H^1(\xi,l), v(l)=0\}.$
Now we are in position to establish the following main result.
\begin{thm}\label{thm1}
For any initial data $(X^0,u_1^0,u_2^0)\in Y$ satisfying the following compatibility conditions
\begin{equation}\label{cond1}
u_1^0(\xi)-\int_0^\xi k_1(\xi,y)u_1^0(y)dy+\vf(\xi)X^0=u_2^0(\xi)+\int_{\xi}^l k_2(\xi,y)u_2^0(y)dy
\end{equation}
and
\begin{align}\label{cond2}
(u_1^0)^{'}(\xi)-k_1(\xi,\xi)u_1^0(\xi)-\int_0^\xi \frac{\p k_1}{\p x}(\xi,y)u_1^0(y)dy+\vf^{'}(\xi)X^0\nonumber\\
 =(u_2^0)^{'}(\xi)-k_2(\xi,\xi)u_2^0(\xi)+\int_{\xi}^l\frac{\p k_2}{\p x}(\xi,y)u_2^0(y)dy,
\end{align}
system (\ref{sys2}) with the feedback law (\ref{fedb}) has a unique  classical solution in
$C\big([0,+\infty[,Y)\cap C^1\big([0,+\infty[,H\big).$
Moreover,  there exists $C > 0$ and $d > 0$   such that the solution satisfies
\begin{equation}\label{expst}
\|\big(X(t),u_1(.,t),u_2(.,t)\big)\|_Y\leq C e^{-d\displaystyle t}\|\big(X^0,u_1^0,u_2^0\big)\|_Y,\quad \forall t\geq 0.
\end{equation}
\end{thm}
\begin{rem}\label{rem1}
 In particular, the compatibility condition (\ref{cond1}) implies that $w_1^0(\xi)= w_2^0(\xi)$. Hence, by property of semi-groups, we get immediately
 $$w_1(t,\xi)=w_2(\xi,t),\quad \forall t\geq 0.$$
\end{rem}
\section{Backstepping Design}
We are now going to find  the gain kernels $ k(x,y)$ and  $k_2(x,y)$ and the gain function $\vf(x)^T.$
From (\ref{sys2}), (\ref{w2}) and (\ref{target}), we get
\begin{eqnarray}\label{w1t}
\frac{\p w_1}{\p t}(x,t)&=&\frac{\p u_1}{\p t}(x,t)-\int_0^x k_1(x,y)\frac{\p u_1}{\p t}(y,t)dy+\vf(x)\dot{X}(t)\nonumber\\
&=& \frac{\p^2 u_1}{\p x^2}(x,t)+\lambda u_1(x,t)-\int_0^x k_1(x,y)\big(\frac{\p^2 u_1}{\p y^2}(y,t)+\lambda u_1(y,t)\big)dy+\vf(x)\dot{X}(t)\nonumber\\
&=& \frac{\p^2 u_1}{\p x^2}(x,t)+\lambda u_1(x,t)-k_1(x,x)\frac{\p u_1}{\p x}(x,t)+k_1(x,0)\frac{\p u_1}{\p x}(0,t)+ \frac{\p k_1}{\p y}(x,x)u_1(x,t)\nonumber\\
&& -\int_0^x \big(\frac{\p^2 k_1}{\p y^2}(x,y)+\lambda k_1(x,y)\big)u_1(y,t)dy+\vf(x)AX(t)+\vf(x)B\frac{\p u_1}{\p x}(0,t),
\end{eqnarray}
\begin{equation}\label{w1x}
\frac{\p w_1}{\p x}(x,t)= \frac{\p u_1}{\p x}(x,t)-k_1(x,x)u_1(x,t)-\int_0^x \frac{\p k_1}{\p x}(x,y)u_1(y,t) dy+\vf^{'}(x)X(t),
\end{equation}
and
\begin{eqnarray}\label{w1xx}
\frac{\p^2 w_1}{\p x^2}(z,t)&=& \frac{\p^2 u_1}{\p x^2}(x,t)-\frac{d}{dx}(k_1(x,x))u_1(x,t)-k_1(x,x)\frac{\p u_1}{\p x}(x,t)
-\frac{\p k_1}{\p x}(x,x)u_1(x,t)\nonumber\\
&&-\int_0^x \frac{\p^2 k_1}{\p x^2}(x,y)u_1(y,t)dy+\vf^{''}(x)X(t).
\end{eqnarray}
Combining (\ref{w1t}) and (\ref{w1xx}) gives
\begin{eqnarray}\label{w1zero}
0 &=&\frac{\p w_1}{\p t}(x,t)-\frac{\p^2 w_1}{\p x^2}(x,t)\nonumber\\
  &=& \big(\lambda +2\frac{d}{d x}(k_1(x,x))\big)u_1(x,t)+\big(k_1(x,0)+\vf(x)B\big)\frac{\p u_1}{\p x}(0,t) -(\vf^{''}(x)-\vf(x)A)X(t)\nonumber\\
 && +\int_0^x\big(\frac{\p^2 k_1}{\p x^2}(x,y) -\frac{\p^2 k_1}{\p y^2}(x,y)-\lambda k_1(x,y)\big)u_1(y,t)dy.
\end{eqnarray}
In the same way, we can just get the following identity for $w_2$
\begin{eqnarray}\label{w2zero}
0 &=&\frac{\p w_2}{\p t}(x,t)-\frac{\p^2 w_2}{\p x^2}(x,t)\nonumber\\
  &=& \big(\lambda +2\frac{d}{d x}(k_2(x,x))\big)u_2(x,t)+k_2(x,l)\frac{\p u_2}{\p x}(l,t)\nonumber\\
 && +\int_x^l\big(\frac{\p^2 k_2}{\p y^2}(x,y) -\frac{\p^2 k_2}{\p x^2}(x,y)+\lambda k_2(x,y)\big)u_2(y,t)dy.\nonumber\\
\end{eqnarray}
Moreover, setting $x=0$ in $w_1(x,t)$, $x=l$ in $w_2(x,t)$ and $\xi$ in both $\frac{\partial w_1}{\partial x}(x,t)$ and $\frac{\partial w_2}{\partial x}(x,t)$,
 and taking Remark \ref{rem1} in mind,
it follows that, if the gain function $\vf(x)^T$ defined in $[0, l]$, the gain kernel $q_1(x,y)$ defined in
\begin{equation}\label{T1}
 T_1=\{(x,y)\mid x\in [0,\xi],y\in [x,\xi]\}
 \end{equation}
and  the gain kernel $q_2(x,y)$ defined in
\begin{equation}\label{T2}
 T_2=\{(x,y)\mid x\in [\xi,l],y\in [x,l]\}
 \end{equation}
 satisfy
\begin{equation}\label{fisys}
\left\lbrace\begin{array}{lll}
\vf^{''}(x)-\vf(x)A=0,\\
\vf(0)=0,\\
\vf^{'}(0)=-K,
 \end{array}\right.
\end{equation}
\begin{equation}\label{k1sys}
\left\lbrace\begin{array}{lll}
\frac{\p^2 k_1}{\p x^2}(x,y)-\frac{\p^2 k_1}{\p y^2}(x,y)=\lambda k_1(x,y),\\
k_1(x,0)= -\vf(x)B,\\
k_1(x,x)=-\frac{\lambda}{2}x,
 \end{array}\right.
\end{equation}
and
\begin{equation}\label{k2sys}
\left\lbrace\begin{array}{lll}
\frac{\p^2 k_2}{\p x^2}(x,y)-\frac{\p^2 k_1}{\p y^2}(x,y)=-\lambda k_2(x,y),\\
k_2(x,l)= 0,\\
k_2(x,x)=\frac{\lambda}{2}(l-x),
 \end{array}\right.
\end{equation}
respectively, then  we obtain  the target system (\ref{target}) for every solution of
the closed loop (\ref{sys2}) with the feedback law (\ref{fedb}).
Obviously, the solution of the linear differential equation (\ref{fisys}) is
\begin{equation}\label{fisol}
\vf(x)=(0,-K)e^{x M}E,
\end{equation}
where $M$ and $E$ are the constant matrices $$M=\begin{pmatrix}
0 & A\\
I_n & 0
\end{pmatrix}, \;\; E=\begin{pmatrix}
I_n \\
0
\end{pmatrix}.$$
According to \cite{Krisbookcourse}, the solution of (\ref{k2sys}) can be done explicitly. For (\ref{k1sys}),
because of the presence of  $\vf(x)$ in the boundary, we can only prove the existence of $k_1(x,y)$.
To study (\ref{k1sys}), we first convert it into an integral equation. For this raison, introducing the change of variables
$$\zeta=x+y\,\,\,\,,\,\,\,\,\eta=x-y,$$
and define $$G(\zeta,\eta):= q(x,y).$$
Then, the function $G$ defined in the triangle $T_0=\Big\{(\zeta,\eta),\eta\in[0,l],\zeta\in[\eta,2l-\eta]\Big\}$, satisfies
\begin{equation}\label{Gsys}
\left\lbrace\begin{array}{lll}
\frac{\p^2 G}{\p\zeta\p\eta}(\zeta,\eta)=\frac{\lambda}{4}G(\zeta,\eta),\,\, (\zeta,\eta) \in T_0,\\
G(\eta,\eta)= -\vf(\eta)B,\,\,\eta \in [0,l],\\
G(\zeta,0)=-\frac{\lambda}{4}\zeta,\,\, \zeta \in [0,2l].
 \end{array}\right.
\end{equation}
Integrating the first equation of (\ref{Gsys}) with respect to $\eta$ from $0$ to $\eta$ and using the boundary condition $G(\zeta,0)$, we get
$$\frac{\p G}{\p\zeta}(\zeta,\eta)=-\frac{\lambda}{4}+\int_0^{\eta}\frac{\lambda}{4}G(\zeta,s)ds.$$
Next, integrating the above identity with respect to $\zeta$ over the interval $[\eta,\zeta]$ and using the boundary condition $G(\eta,\eta)$,
it follows
\begin{equation}
G(\zeta,\eta)= -\vf(\eta)B-\frac{\lambda}{4}(\zeta-\eta)+\int_{\eta}^{\zeta}\int_0^{\eta}\frac{\lambda}{4}G(t,s)ds dt,\forall n\in\N.
\end{equation}
To achieve the existence of $G$ in $T_0$, we use the method of successive approximations. To this end, let us set
\begin{eqnarray}
G^0(\zeta,\eta)&=&0,\\
G^{n+1}(\zeta,\eta)&=&-\vf(\eta)B-\frac{\lambda}{4}(\zeta-\eta)+\int_{\eta}^{\zeta}\int_0^{\eta}\frac{\lambda}{4}G^n(t,s)ds dt, \forall n\in\N,
\end{eqnarray}
and denote the difference between two consecutive terms by
\begin{equation}\label{deltagn}
\Delta G^n(\zeta,\eta)=G^{n+1}(\zeta,\eta)-G^n(\zeta,\eta).
\end{equation}
Then,
\begin{equation}
\Delta G^{n+1}(\zeta,\eta)=\frac{\lambda}{4}\int_{\eta}^{\zeta}\int_0^{\eta}\Delta G^n(t,s)ds dt.
\end{equation}
We have
$$\Delta G^0(\zeta,\eta)=-\frac{\lambda}{4}(\zeta-\eta)-\vf(\eta)B.$$
Since $\vf$ is continuous on $[0,l]$, there exists  $\mu>0$ such that $|\vf(\eta)B|\leq \mu, \forall \eta\in [0,l].$
Thus, $$|\Delta G^0(\zeta,\eta)|\leq \frac{\lambda}{4}(\zeta-\eta)+\mu.$$
Using the fact that $0\leq\eta\leq\zeta$, by an immediate mathematical induction, it can be shown that
\begin{equation}\label{estdeltaG}
|\Delta G^n(\zeta,\eta)|\leq (\frac{\lambda}{4})^{n+1}\frac{(\zeta-\eta)\zeta^n\eta^n}{n!(n+1)!}+\mu (\frac{\lambda}{4})^n\frac{\zeta^n\eta^n}{(n!)^2}.
\end{equation}
It then follows from the Weierstrass M-test that the series
$$G(\zeta,\eta)=\lim_{n\rightarrow \infty}G^n(\zeta,\eta)=\sum_{n=0}^{\infty}\Delta G^n(\zeta,\eta)$$
converges absolutely and uniformly in $T_0$. Having proved the existence of $G(\zeta,\eta)$ in $T_0$, that of $k_1(x,y)$ in $T$ follows immediately.
For the kernel\,\, $k_2(x,y)$, without loss of generality, we suppose that $l-\xi\leq\xi$. Consider the change of coordinates
$$s=l-x\,\,\,,\,\,\,t=l-y,$$ that maps the triangle $T_2$ into the triangle $T_1$, and
define the function $h(s,t)=-k_2(x,y)$. Then, we get for the function \, $h(s,t)$ the following PDE
\begin{equation}\label{hsys}
\left\lbrace\begin{array}{lll}
\frac{\p^2 h}{\p s^2}(s,t)-\frac{\p^2 h}{\p t^2}(s,t)=\lambda h(s,t),\\
h(s,0)= 0,\\
h(s,s)=\frac{\lambda}{2}s.
 \end{array}\right.
\end{equation}
According to \cite{Krisbookcourse}, the solution of (\ref{hsys}) in $T_1$ is
\begin{eqnarray}
h(s,t)&=&-\lambda\frac{I_1\Big(\sqrt{\lambda(s^2-t^2)}\Big)}{\sqrt{\lambda(s^2-t^2)}},\,\, if\,\, s\neq t,\nonumber\\
h(s,s)&=&\frac{\lambda}{2}s,
\end{eqnarray}
where $I_1$  is the first order modified Bessel function of the first kind.
Consequently, we obtained the solution of (\ref{k2sys}) in $T_2$ as follows
\begin{eqnarray}
k_2(x,y)&=&\lambda\frac{I_1\Big(\sqrt{\lambda ((l-x)^2-(l-y)^2})\Big)}{\sqrt{\lambda ((l-x)^2-(l-y)^2)}},\,\, if\,\, x\neq y. \nonumber\\
k_2(x,x)&=&\frac{\lambda}{2}(l-x).
\end{eqnarray}
Let's move to the proof of the main Theorem \ref{thm1}.

\section{Proof of Theorem {\ref{thm1}}}
\subsection{Well posedness of the initial plant (\ref{sys2})}

Let us define the map
\begin{eqnarray}
  \Omega\,\, :\,\, Y &\rightarrow & Y \nonumber \\
 (X,u_1,u_2)&\mapsto & (X,w_1,w_2),
\end{eqnarray}
where $w_1$ and $w_2$ satisfy (\ref{w1}) and (\ref{w2}), respectively. This linear map is well defined and bounded. Hence, there exists
a positive constant $c_1$ such that
\begin{equation}\label{omegacont}
\|\Omega(X,u_1,u_2)\|_Y\leq c_1 \|(X,u_1,u_2)\|_Y,\,\,\forall (X,u_1,u_2)\in Y.
\end{equation}
It is a straightforward that the transformation $\Omega$ is an isomorphism, and
\begin{eqnarray}
\Omega^{-1} &:& Y\rightarrow Y\\
  \nonumber&& (X,w_1,w_2)\mapsto (X,u_1,u_2)
\end{eqnarray}
has the following form
 \begin{align}
X(t)&=X(t)\\
u_1(x,t)&= w_1(x,t)-\int_0^x w_1(y,t)g_1(x,y)dy+\psi(x)X(t),\quad x\in [0,\xi], \quad t\geq 0,\label{u1}\\
u_2(x,t)&= w_2(x,t)+\int_x^l w_2(y,t)g_2(x,y)dy, \quad x\in [\xi,l], \quad t\geq 0.\label{u2}
\end{align}
As is done in the study of the direct transformation $\Omega$ and in the same way, on can prove the existence of the gain kernel
 $g_1(x,y)$ and compute explicitly  the gain kernel  $g_2(x,y)$ and the gain function $\psi(x)$.
 Therefore, there exists a positive constant $c_2$ such that
 \begin{equation}\label{omega-1cont}
\|\Omega^{-1}(X,w_1,w_2)\|_Y\leq c_2 \|(X,w_1,w_2)\|_Y,\,\,\forall (X,w_1,w_2)\in Y.
 \end{equation}
It is well known that if the initial condition $(w_1^0,w_2^0)$ of the target system (\ref{target}) belongs to the subspace
$$D=\left\{(w_1,w_2)\in H^2(0,\xi)\cap H_L^1(0,\xi)\times H^2(\xi,l)\cap H_R^1(\xi,l)\mid w_1(\xi)=w_2(\xi),\frac{\p w_1}{\p x}(\xi)=
\frac{\p w_2}{\p x}(\xi)\right\},$$
the function $w$ defined as
$$w(x,t)=
\left\lbrace\begin{array}{lll}
w_1(x,t), \quad x\in [0,\xi], \quad t\geq 0,\\
w_2(x,t), \quad x\in [\xi,l], \quad t\geq 0.
\end{array}\right.$$
belongs to $H^2(0,l)\cap H_0^1(0,l)$ and satisfying the following boundary problem
\begin{equation}\label{wsys}
\left\lbrace\begin{array}{lll}
\frac{\partial w}{\partial t}(x,t) = \frac{\partial^2 w}{\partial x^2}(x,t),\quad t>0, \quad x\in(0,l),\\
 w(0,t)=w(l,t)=0, \quad t>0,\\
 w(x,0)=w^0(x)=\mathds{1}_{(0,\xi)}w_1^0(x)+\mathds{1}_{(\xi,l)}w_2^0(x),\quad x\in(0,l),
\end{array}\right.
\end{equation}
The existence , uniqueness and regularity of the solution of (\ref{wsys}) follow by standard
arguments of semi-groups theory. Furthermore, by the method of separation of variables, the system (\ref{wsys}) is not only well posed but its solution
 is explicitly done as follows
\begin{equation}\label{wexpl}
w(x,t)=\frac{2}{l}\sum_{k=1}^{+\infty}e^{-\frac{k^2\pi^2}{l^2}t}sin\Big(\frac{k\pi}{l}x\Big)\int_0^l w^0(s)sin\Big(\frac{k\pi}{l}s\Big)ds.
\end{equation}
Thus, By Duhamel's formula, the solution of the ODE in the target system (\ref{target}) yields
\begin{equation}
X(t)= e^{t(A+BK)}X_0+\int_0^t e^{(t-\tau)(A+BK)}B\frac{\p w_1}{\p x}(0,\tau)d\tau.
\end{equation}
Since the isomorphism $\Omega$ transforms system (\ref{sys2}) to system (\ref{target}), it follows that system (\ref{sys2}) with the feedback
law (\ref{fedb}) is well posed. Hence, the regularity of the solution given by Theorem \ref{thm1} holds true.
\subsection{Exponential stability}

Consider the Lyapunov function candidate
\begin{equation}\label{Vly}
V(t)= X(t)^TPX(t)+\frac{a}{2}\|w(.,t)\|^2+\frac{b}{2}\|\frac{\p w}{\ x}(.,t)\|^2,
\end{equation}
where $a>0$ and $b>0$ are two constants which we will specify later and the positive definite matrix $P=P^T>0$ is the solution of the Lyapunov equation
$$ P(A+BK)+(A+BK)^T P = -Q,$$
for some positive definite matrix $Q=Q^T>0$.  From (\ref{Vly}), it can be obtained that for all $t\geq 0$,
\begin{equation}\label{Venc}
 \alpha_1\|(X(t), w(.,t))\|_Z^2\leq V(t)\leq \alpha_2\|(X(t), w(.,t))\|_Z^2,
\end{equation}
where
 $$\alpha_1=\min\big(\la_{min}(P),\frac{a}{2},\frac{b}{2})\,\,\,\,\, , \,\,\,\,\,\alpha_2=\max\big(\la_{max}(P),\frac{a}{2},\frac{b}{2}\big)$$
and  $Z=\R^n\times H^1(0,1)$ equipped with its norm $\|(X, w)\|_Z^2=|X|^2+\| w\|_{H^1(0,1)}^2$.
The derivative of $V$ along the solutions of (\ref{wsys}) is given by
$$\dot{V}(t)=-X^T(t)QX(t)+2X(t)^T PB\frac{\p w}{\p x}(0,t)-a \|\frac{\p w}{\p x}(.,t)\|^2-b\|\frac{\p^2 w}{\p^2 x}(.,t)\|^2.$$
By Young's inequality, we get
$$2X(t)^T PB \frac{\p w}{\p x}(0,t)\leq \frac{\la_{min}(Q)}{2}|X(t)|^2+\frac{2}{\la_{min}(Q)}|PB|^2 (\frac{\p w}{\p x}(0,t))^2,$$
and by Agmon's inequality \cite{{Tang11SCL}}, it can be proved that the following inequality holds:
$$-\|\frac{\p^2 w}{\p^2 x}(.,t)\|^2\leq \frac{1+l}{l}\|\frac{\p w}{\p x}(.,t)\|^2-(\frac{\p w}{\p x}(0,t))^2.$$
Thus, $$ \dot{V}(t)\leq -\frac{\la_{min}(Q)}{2}|X(t)|^2-\Big(\frac{a}{2}-\frac{b(1+l)}{l}\Big)\|\frac{\p w}{\p x}(.,t)\|^2
-\frac{a}{2}\|\frac{\p w}{\p x}(.,t)\|^2-\Big(b-\frac{2|PB|^2}{\la_{min}(Q)}\Big)(\frac{\p w}{\p x}(0,t))^2.$$
Since $w(0,t)=0$, by Poincar\'e inequality, we get
$$ \|w(.,t)\|^2\leq 4l^2\|\frac{\p w}{\ x}(.,t)\|^2.$$
Now, if we choose
$$b>\frac{2|PB|^2}{\la_{min}(Q)}\quad and \quad a>\frac{2b(1+l)}{l}+2,$$
from estimation above, it yields
\begin{eqnarray}
\dot{V}(t)&\leq&-\frac{\la_{min}(Q)}{2}|X(t)|^2-\frac{a}{8l^2}\|w(.,t)\|^2-\|\frac{\p w}{\p x}(.,t)\|^2\nonumber\\
&\leq&-\delta V(t),
\end{eqnarray}
where $$\delta= min\Big(\frac{\la_{min}(Q)}{2\la_{max}(P)},\frac{1}{4l^2},\frac{2}{b}\Big).$$
Therefore,
\begin{equation}\label{V0}
V(t)\leq V(0)e^{-\delta \displaystyle t}, \forall t\geq 0.
\end{equation}
 Consequently, From (\ref{Venc}) and (\ref{V0}), the following estimation
\begin{equation}\label{wexps}
\|(X(t), w(.,t))\|_Z^2\leq \alpha e^{-\delta \displaystyle t}\|(X^0, w^0\|_Z^2
\end{equation}
 holds for all $t\geq0$, where $\alpha=\frac{\alpha_1}{\alpha_2}$.\\
 Ultimately, since $\|(X,w_1,w_2)\|_Y= \|(X, w)\|_Z$, from (\ref{omegacont}), (\ref{omega-1cont}) and (\ref{wexps}), it follows
\begin{equation}
\|(X(t),u_1(.,t),u_2(.,t)\|_Y\leq C e^{-d\displaystyle t}\|\big(X^0,u_1^0,u_2^0\big)\|_Y,\quad \forall t\geq0,
\end{equation}
 where $C=c_1 c_2\sqrt{\alpha}$ and $d=\frac{\delta}{2}.$
Thus, the proof of Theorem \ref{thm1} is complete.

\section{Conclusion and future work}
We have considered the exponential stabilization for a cascaded ODE-reaction-diffusion system (\ref{sys1}) coupling at the left boundary
by Neumann connection and a control acting at an internal point of the PDE subdomain. By the backstepping method, we have construct a feedback
law to achieve the result. However for a general coupled system where the cascade ODE-PDE system is coupled at an intermediate point $x_0\in(0,l)$ as follows
\begin{equation}\label{sysfin}
\left\lbrace\begin{array}{lll}
\dot{X}(t)=AX(t)+B\frac{\p u}{\p x}(x_0,t),\quad t>0,\\
\frac{\p u}{\p t}(x,t) = \frac{\partial^2 u}{\partial x^2}(x,t)+\lambda u(x,t) ,\quad t>0, \quad x\in(0,\xi)\cup (\xi,l),\\
u(\xi^-,t) = u(\xi^+,t) ,\quad t>0,\\
\frac{\partial u}{\partial x}(\xi^-,t)-\frac{\partial u}{\partial x}(\xi^+,t)=U(t) ,\quad t>0,\\
X(0)= X^0,\\
u(0,t)  = u(l,t)= 0,\quad t>0,\\
u(x,0)=u^0(x) ,\quad x\in (0,l),
\end{array}\right.
\end{equation}
the stabilization controller design is not obvious, which needs to be considered in the future.


\begin{thebibliography}{999}

\bibitem{Am01hetu} K. Ammari, A. Henrot and M. Tucsnak, Asymptotic behaviour of the solutions and optimal location of the actuator for the pointwise stabilization of a string, Asymptotic Anal., 28 (2001), 215-2ֲ40.

\bibitem{Am00tu} K. Ammari and M. Tucsnak, Stabilization of Bernoulli-Euler beams by means of a pointwise feedback force, SIAM J. Control Optim., 39 (2000), 1160-1181.

\bibitem{Coron98} J.-M. Coron and B. D'Andr\'ea Novel, Stabilization of a rotating body beam without damping, IEEE Transaction on Automatic Control,
43 (1998), 608-618.

\bibitem{Fa17} F. Hassine, Rapid exponential stabilization of 1-d transmission wave equation with in-domain anti-damping, Asian Journal of Control 19 (2017), 1-11.

\bibitem{Kris09heat} M. Krstic, Compensating actuator and sensor dynamics governed by diffusion PDEs, Systems $\&$ Control Letters, 58 (2009), 372-377.

\bibitem{Krisbook} M. Krstic, Delay Compensation for Nonlinear, Adaptive, and PDE Systems, Birkhauser, (2009).

\bibitem{Kris10JFI} M. Krstic and G.-A. Susto, Control of PDE-ODE cascades with Neumann interconnections, Journal of Franklin Institute,
347 (2010), 284-314.

\bibitem{Kris95} M. Krstic, I. Kanellakopoulos and P. Kokotovic, Nonlinear and Adaptive Control Design, John Wiley and Sons, (1995).


\bibitem{Kris00} M. Krstic and W. Liu, Backstepping boundary control of Burgers equation with actuator dynamics, Systems $\&$ Control Letters,
 41 (2000), 291-303.

\bibitem{Kris04} M. Krstic and A. Smyshlyaev, Closed-form boundary state feedbacks for a class of $1-d$ partial
integro-differential equations, IEEE Transaction on Automatic Control, 49 (2004), 2185-2202.

\bibitem{Krisbookcourse} M. Krstic and A. Smyshlyaev, Boundary Control of PDEs: A Course on Backstepping Designs, SIAM, (2008).

\bibitem{Liu03} W. Liu, Boundary feedback stabilization of an unstable heat equation, SIAM Journal Control Optimization, 42 (2003), 1033-1043.

\bibitem{Tang11} S. Tang and C. Xie, Stabilization for acoupled PDE-ODE control system, Journal of the Franklin Institute, 348 (2011), 2142-2155.

\bibitem{Zhoua12} S. Tang and Z. Zhoua, Boundary stabilization of a coupled wave-ODE system with internal anti-damping, International Journal of Control, 85 (2012), 1683-1693.

\bibitem{Tang11SCL} S. Tang and C. Xie,  State and output feedback boundary control for a coupled PDE-ODE system.
Systems $\&$ Control Letters, 60 (2011), 540-545.

\bibitem{Woi13} S. Wang, and F. Woittennek, Backstepping-method for parabolic systems with in-domain actuation, IFAC, Proc, 1 (2013), 43-48.

\end{thebibliography}
\end{document}